\documentclass[12pt, notitlepage]{article} 
\usepackage{amsmath, amsfonts, amscd, latexsym, amsthm, amssymb, graphicx}
\newtheorem{theor}{\hspace{1cm}{\sc Theorem}}[section]
\newtheorem{utver}[theor]{\hspace{1cm}{\sc Proposition}}
\newtheorem{sledst}[theor]{\hspace{1cm}{\sc Corollary}}
\newtheorem{lemma}[theor]{\hspace{1cm}{\sc Lemma}}
\theoremstyle{definition}
\newtheorem{defin}[theor]{\hspace{1cm}{\sc Definition}}
\newtheorem*{rem}{\hspace{1cm}{\sc Remark}}
\newtheorem*{exa}{\hspace{1cm}{\sc Example}}

\newtheorem*{conjec}{\hspace{1cm}{\sc Conjecture}}

\newcommand{\Vol}{\mathop{\rm Vol}\nolimits}

\newcommand{\codim}{\mathop{\rm codim}\nolimits}

\def\R{\mathbb R}

\def\Z{\mathbb Z}
\def\Q{\mathbb Q}
\def\C{\mathbb C}
\def\CC{({\mathbb C}\setminus\{0\})}
\def\T{\mathbb T}

\begin{document}
 
\begin{center}
{\Large \textbf{Tropical varieties with polynomial weights\\ and
corner loci of piecewise polynomials.}}
\end{center}
 
\begin{center}
A. Esterov\footnote{\noindent Partially supported by
RFBR-10-01-00678 grant. E-mail address: \texttt{esterov@mccme.ru}}
\end{center}
\vspace{-1cm}
\begin{flushright} \textit{To S. M. Gusein-Zade on \\ the occasion of his
60th birthday}
\end{flushright}
 
\vspace{-1cm}
\section{Introduction.}
 
Counting Euler characteristics of the discriminant of the
quadratic equation in terms of Newton polytopes in two different
ways, G. Gusev (\cite{gphd}) found an unexpected relation for
mixed volumes of two polytopes $S_1$ and $S_2\subset\R^n$ and the
convex hull $S$ of their union. For instance, assuming $n=2$ and
denoting the mixed area of polygons $P$ and $Q$ by
$\Vol(P,Q)=\Vol(P+Q)-\Vol(P)-\Vol(Q)$, this relation specializes
to
$$\Vol(S,S)-\Vol(S,S_1)-\Vol(S,S_2)+\Vol(S_1,S_2)=0.$$ We call it
unexpected because it is not a priori invariant under parallel
translations of $S_1$.

We give an elementary proof and a multidimensional generalization
of this equality as requested in \cite{gphd} (see Corollary \ref{sled1} below),
deducing it from the following fact (Theorem \ref{prodsup0}):
the mixed volume of
polytopes depends only on the product of their support functions, rather than on individual support functions.
We give a new elementary formula for this dependence (Proposition
 \ref{prodsup01}) and represent 
it as a specialization of the isomorphism
between two well known combinatorial models of the cohomology of
toric varieties. 
 
Although this construction makes sense 
for arbitrary polytopes, it so far has been established only for polytopes
with rational vertices (partially due to the lack of combinatorial
tools capable of substituting for toric geometry, when vertices are not rational),
 see e. g. \cite{katz1}.
To fill this gap, we introduce tropical varieties with polynomial weights,
i.e. fans with somehow balanced polynomial functions on their cones (see Definition \ref{deftrop}). 
 
 This notion 
interpolates between the notions of conventional tropical varieties and continuous
piecewise polynomial functions. It allows us to establish the aforementioned results 
 for non-rational polytopes.
 
We also discuss possible applications of polynomially weighted tropical varieties 
to tropical intersection theory. Namely, we notice that the intersection theory
 on a smooth tropical fan, recently constructed in \cite{all1}, \cite{fr}, \cite{sh1},
 can be seen as the restriction of the intersection theory on the ambient vector space (see Theorem \ref{smoothhh}).
 Polynomially weighted tropical varieties allow to conjecture a generalization of
 this fact to non-smooth tropical varieties.
 
 The four preceding paragraphs describe the contents of the four sections of the paper.
 
\textbf{Gusev's equality.} To simplify notation, we denote the
mixed volume of polytopes $A_1,\ldots,A_k$ in $\R^k$ by the
monomial $A_1\cdot\ldots\cdot A_k$ (this mixed volume is by
definition the coefficient of the monomial $x_1\ldots x_k$ in the
polynomial $\Vol(x_1A_1+\ldots+x_kA_k)$ of variables
$x_1,\ldots,x_k$). In the same way, for a homogeneous polynomial
$P(x_1,\ldots,x_m)=\sum c_{a_1,\ldots,a_m}x_1^{a_1}\ldots
x_m^{a_m}$ of degree $k$, we define $P(A_1,\ldots,A_m)$ as $\sum
c_{a_1,\ldots,a_m}A_1^{a_1}\ldots A_m^{a_m}$.
\begin{theor}[\cite{gphd}] \label{gusev1}
For any two polytopes $S_1$ and $S_2\subset\R^n$ and the convex
hull $S$ of their union, we have
$(2^n-2)S^n=\sum_{i=1}^{n-1}2^i\bigl(S_1^{n-i}S^i+S^{n-i}S_2^i-S_1^{n-i}S_2^i\bigr)$.
\end{theor}
We deduce this from the following fact. Denote the support
function of a polytope $A\subset\R^n$ by
$A(\cdot):(\R^n)^*\to\R$, so that $A(v)=\max_{a\in A} v\cdot a$.
\begin{theor}\label{prodsup0} There exists a linear function $D$ on the space of
conewise-polynomial functions on $(\R^n)^*$, such that
$$D\Bigl(A_1(\cdot)\ldots A_n(\cdot)\Bigr)=A_1\ldots A_n$$ for
every collection of polytopes $A_1,\ldots,A_n$ in $\R^n$.
\end{theor}
Recall that a function $f$  on $\R^m$ is said to be conewise-polynomial, if
it is polynomial on every piece of a finite subdivision of $\R^m$
into polyhedral cones with vertices at $0$. \par
 
 Note that the existence of a function $D$ (aside from its
linearity) is not obvious a priori, because the collection of
polytopes is not uniquely determined by the product of their
support functions: the two different pairs of polygons on the
following picture have the same product of support functions
(and, thus, the same mixed volume, which is equal to 4).
\begin{center} \noindent\includegraphics[width=12cm]{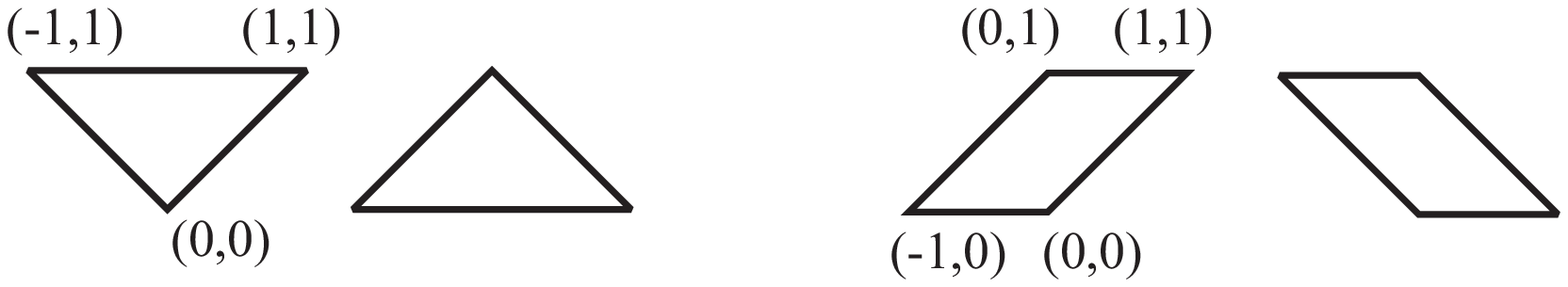}
 
\nopagebreak\small{Picture 1.}
\end{center}
 
Also note that the function $D$ is not monotonous: if $A$, $B$ and $C$ are the segments in the plane from the origin to the points $(1,0)$, $(0,1)$ and $(1,1)$ respectively, then $A(\cdot)B(\cdot)<C(\cdot)^2$, although $A\cdot B=1>C\cdot C=0$. \par 
 
 For rational polytopes, Theorem \ref{prodsup0} is a special case of the isomorphism between two well known models of cohomology of toric varieties, as explained at the end of this section. It also follows from a stronger fact about the product of support functions of rational polytopes: Theorem 5.1 in \cite{katz1}. This cannot be extended to non-rational cones and polytopes by continuity arguments, and we deduce Theorem \ref{prodsup0} in full generality from Proposition \ref{prodsup01} below, suggesting an explicit formula for $D$. \par

 Note that results of \cite{katz1} remain valid for non-rational polytopes as well; to prove them in the non-rational setting,
one should replace the reference to Brion's formula in \cite{katz1} with the reference to the combinatorial
Riemann-Roch formula of \cite{khp} (i.e. to replace
 summation over lattice points of a polyhedron with integration over a polyhedron). However, this is beyond the scope of our paper.
 \par

For an (ordered) basis $v_1,\ldots,v_n$ in $\R^n$, denote the cone generated by $v_1,\ldots,v_n$ by $\langle v_1,\ldots,v_n\rangle$,
and denote the Gram-Schmidt orthogonalization of $v_1,\ldots,v_n$ by $v_1^\bot,\ldots,v_n^\bot$ (so that $v_1^\bot,\ldots,v_n^\bot$
is orthonormal, $v_1^\bot,\ldots,v_i^\bot$ generate the same subspace as $v_1,\ldots,v_i$ do for $i=1,\ldots,n$, and $v_i\cdot v_i^\bot>0$).
For a continuous conewise polynomial function $f:\R^n\to\R$, consider a simple complete fan $\Gamma$, on whose cones $C\in\Gamma$
the function $f$ coincide with polynomials $f_C$. Then Theorem \ref{prodsup0} can be formulated as follows.
\begin{utver}\label{prodsup01} Define $D(f)$ as
$$\frac{1}{n!}\sum_{\langle v_1,\ldots,v_n\rangle\in\Gamma}\frac{\partial^n f_{\langle v_1,\ldots,v_n\rangle}}
{\partial v_1^\bot\ldots\partial v_n^\bot},\eqno{(*)}$$ where the sum is taken over all ordered bases of unit vectors $v_1,\ldots,v_n$,
generating cones from $\Gamma$. Then $D(f)$ does not depend on the choice of the fan $\Gamma$, linearly depends on $f$, and $D\Bigl(A_1(\cdot)\ldots A_n(\cdot)\Bigr)$ equals the mixed volume of polytopes $A_1,\ldots,A_n$.
\end{utver}
 It is an elementary rephrasing of Theorem \ref{mainth2}, whose formulation and proof make use of more general machinery, developed in subsequent sections. For the convenience of the reader, we also outline an elementary proof of Proposition \ref{prodsup01} here. Also note that another explicit formula for $D$ is given in \cite{maz1}. \par
 
\textsc{Sketch of the proof.} Independence of subdivisions of $\Gamma$ and linearity follow by definition, so we only need to check that $D\Bigl(A_1(\cdot)\ldots A_n(\cdot)\Bigr)=A_1\cdot\ldots\cdot A_n$.
Moreover, since $D\Bigl(A_1(\cdot)\ldots A_n(\cdot)\Bigr)$ is a multilinear function of $A_1,\ldots,A_n$,
it is enough to check the equality for $A_1=\ldots=A_n$, i.e. to prove that $D\Bigl(A^n(\cdot)\Bigr)$ equals the volume of the polytope $A$.\par

 As a simplicial chain, the polytope $A$ can be represented as a linear combination of simplices with coefficients $\pm 1$, whose volumes coincide with the terms of the sum $(*)$ up to the signs of their coefficients in the linear combination. This fact implies the desired equality, so it is enough to construct requested simplices. We illustrate this construction, assuming for simplicity that the orthogonal complement to the affine span of every (relatively open) face $B\subset A$ intersects $B$. \par
 
Let $\Gamma_1\subset\R^n$ be the set of all external normal vectors to the faces of $A$ of positive dimension,
and let $\Gamma_2$ be the union of all rays from the origin, passing through the points of faces of $A$ of codimension greater than 1. Then $\Gamma_1\cup\Gamma_2$ subdivides $A$ into $n$-dimensional simplices that are in one to one correspondence with the terms of the sum $(*)$,
 and these terms are equal to the volumes of the corresponding simplices by the subsequent Lemma \ref{tetr1}.
$\quad\Box$\par
 
 \begin{lemma}\label{tetr1} If a simplex in $\R^n$ has $n$ mutually perpendicular edges, then its volume equals the product of their lengths times $1/n!$.
 \end{lemma}
 
\begin{sledst}\label{sled1} For any polytopes $B_1,\ldots,B_n$ in $\R^n$ and the convex
hull $B$ of their union, we have $(B-B_1)\ldots(B-B_n)=0$.
\end{sledst}
\textsc{Proof.} Since $B(v)=\max_i B_i(v)$ for every
$v\in(\R^n)^*$, we have $\Bigl(B(v)-B_1(v)\Bigr)\ldots$ $
\Bigl(B(v)-B_n(v)\Bigr)=0$, and the desired equality follows by
Theorem \ref{prodsup0}. $\quad\Box$
 
\textsc{Proof of Theorem \ref{gusev1}.} Sum up the equality
$2^i(S^{n-i}-S_1^{n-i})(S^i-S_2^i)=0$ (which is a special case of
 Corollary \ref{sled1}) over $i=1,\ldots,n-1$. $\quad\Box$
 
We now show that Theorem \ref{prodsup0} is a special case of the isomorphism
between two well known models for cohomology of toric varieties. 
 
\textbf{Cohomology ring of toric varieties and its Brion-Stanley
description.}
 
The set of all complete rational fans in $\R^n$ admits the
following partial order relation: $\Gamma_1\leqslant\Gamma_2$ if
every cone of the fan $\Gamma_2$ is contained in a cone of the
fan $\Gamma_1$. Denoting the toric variety of a fan $\Gamma$ by
$\T^{\Gamma}$, the natural mapping
$\T^{\Gamma_2}\to\T^{\Gamma_1}$ induces a homomorphism of
cohomology rings
$h_{\Gamma_1,\Gamma_2}:H^{\cdot}(\T^{\Gamma_1})\to
H^{\cdot}(\T^{\Gamma_2})$. The direct system of these rings and
homomorphisms
gives rise to the direct 
limit $$\mathcal{H}=\lim_{\to}H^{\cdot}(\T^{\Gamma}).$$ Note that we get the same ring
$\mathcal{H}$, independently of which version of cohomology theory
we consider (e.g. singular cohomology, Chow cohomology or intersection
cohomology; see e.g. \cite{payne} for a good overview of this kind of results).
There are two well known ways to describe this ring combinatorially.
 
Brion's description of Chow rings \cite{brion} and Stanley's
description \cite{stanley} of intersection cohomology of toric
varieties lead to the following one. Let $\mathcal{P}_{\Q}$ be the
ring of continuous piecewise-polynomial functions on $\R^n$,
whose domains of polynomiality are rational convex polyhedral
cones with the vertex $0$.
Denote its ideal, generated by
linear functions, by $\mathcal{L}_{\Q}$. Then
$\mathcal{H}=\mathcal{P}_{\Q}/\mathcal{L}_{\Q}$.
 
\textbf{Fulton-Kazarnovskii-McMullen-Sturmfels description.}
 
One more combinatorial model for the cohomology ring
$\mathcal{H}$ is given independently by many authors, and is
known as McMullen's polytope weights \cite{mcm}, Fulton--Sturmfels
Minkowski weights \cite{fs}, and Kazarnovskii's c-fans \cite{kaz}.
A $k$-dimensional \textit{weighted piecewise-linear set} is a pair
$(P,p)$, where the \textit{support set} $P\subset\R^n$ is a union
of finitely many rational $k$-dimensional polyhedra (closed and
not necessary bounded), and the \textit{weight} $p:P\to\R$ is a
locally constant function on the set of smooth points of $P$. It is said
to be \textit{homogeneous}, if $P$ is a union of polyhedral cones
with the vertex $0$. For a smooth point $x$ of $P$, let
$N_xP\subset\R^n$ be the codimension $k$ subspace, orthogonal to
the tangent space of $P$ at $x$. The \textit{tropical
intersection number} $\circ_i(P_i,p_i)$ of transversal weighted
piecewise-linear sets $(P_i,p_i)$ with $\sum_i\codim P_i=n$ is
the sum of the products
$\Bigl|\Z^n/\bigoplus_i(\Z^n\cap N_xP_i)\Bigr|\cdot\prod_ip_i(x)$
over all points $x\in\cap_i P_i$ (transversality means that all
$P_i$ are smooth at every point of their intersection, and the
tangent planes are transversal).
 
A weighted piecewise-linear set $(P,p)$ is called a
\textit{tropical variety}, if, for every rational subspace
$L\in\R^n$ of the complementary dimension, the tropical
intersection number $(P,p)\circ(L+x,1)$ does not depend on the
point $x\in\R^n$ (note that the intersection number makes sense
for almost all $x$). Arbitrary tropical varieties $(P_i,p_i)$
with $\sum_i\codim P_i=n$ in $\R^n$ intersect transversally when
shifted by generic vectors $x_i\in\R^n$, and this intersection
number $\circ_i\Bigl(P_i+x_i,\, p_i(\cdot-x_i)\Bigr)$ does not
depend on the choice of $x_i$. This allows to call it the
intersection number of the varieties $(P_i,p_i)$ and to denote it
by $\circ_i(P_i,p_i)$. See, for example, the two ways to count the
intersection number of a pair of tropical curves on the right of
the following picture; both ways lead to the same answer 4.
\begin{center} \noindent\includegraphics[width=15cm]{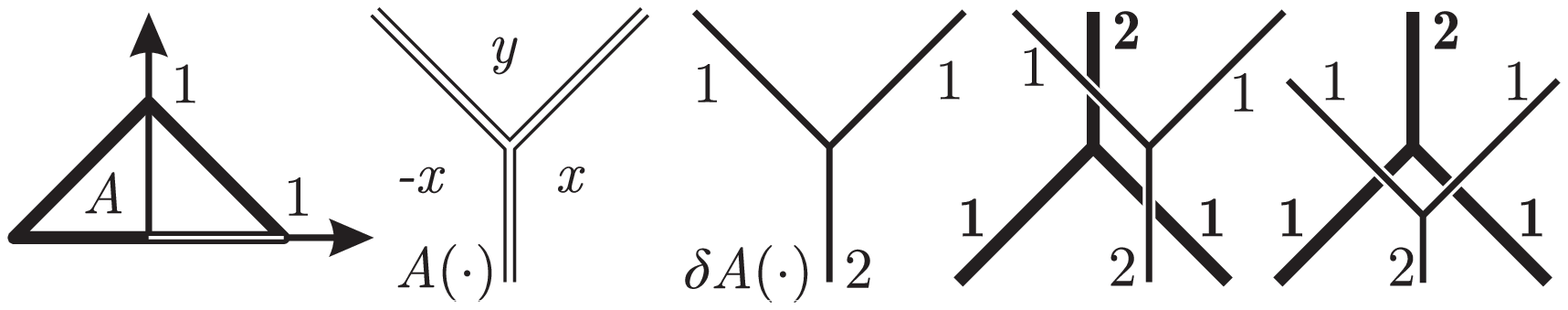}
 
\small{Picture 2.}
\end{center}
The \textit{product} $(R,r)$ of tropical varieties $(P,p)$ and
$(Q,q)$ is uniquely characterized by the equality of the
intersection numbers $(R,r)\circ(S,s)=(P,p)\circ(Q,q)\circ(S,s)$
for every tropical variety $(S,s)$ of the complementary dimension
(the existence of such $(R,r)$ is not clear, see a more
constructive definition in Section 2). In particular, if $(P,p)$
and $(Q,q)$ are homogeneous tropical varieties of complementary
dimension, then their product is the 0-dimensional tropical
variety $\Bigl(\{0\},(P,p)\circ(Q,q)\Bigr)$. With respect to this
multiplication, the natural addition $(P,p)+(Q,q)=(P\cup Q,p+q)$,
and the equivalence relation $(P,0)=(\varnothing,0)$ for every set
$P$, homogeneous tropical varieties form a ring
$\mathcal{C}_{\Q}$, and we have $\mathcal{C}_{\Q}=\mathcal{H}$.
 
\textbf{The isomorphism.}
 
The isomorphisms
$\mathcal{P}_{\Q}/\mathcal{L}_{\Q}=\mathcal{H}=\mathcal{C}_{\Q}$
induce the isomorphism
$I_{\Q}:\mathcal{P}_{\Q}/\mathcal{L}_{\Q}\to\mathcal{C}_{\Q}$ of
the two combinatorial models for cohomology of toric varieties.
There is one more well known combinatorial model for $\mathcal{H}$
by Khovanskii and Pukhlikov, whose isomorphism with $\mathcal{C}_{\Q}$
is combinatorially described in \cite{kk1}, but we do not need this
construction here.
 
Explicit combinatorial constructions for the isomorphism $I_{\Q}$
are given in \cite{katz1} and \cite{maz1}.
Its degree 1 component, sending conewise linear functions to
homogeneous tropical hypersurfaces, is much simpler and admits
the following well known description.
\begin{defin} \label{corner0} Assume that a continuous conewise linear function $L:\R^n\to\R$
equals linear functions $L_+$ and $L_-$ on complementary
half-spaces $H_+$ and $H_-$, separated by a rational hyperplane
$P$ (such a function is called a \textit{book}). Choose a vector
$v\in H_+$ that generates the 1-dimensional lattice $\Z^n/P$, and
define the (constant) function
$$p(x)=\partial_v L_+(x)-\partial_v L_-(x) \mbox{ for every }
x\in P.$$ The \textit{corner locus} of $L$ is defined as the pair
$(P,p)$ for $p\ne 0$ and $(\varnothing,0)$ otherwise (i.e. for
linear $L$). It does not depend on the choice of $v$ and is
denoted by $\delta L$. For an arbitrary continuous piecewise
linear function $L$, whose domains of linearity are rational
polyhedra, its corner locus is the weighted piecewise-linear set
$\delta L$, such that whenever $L$ equals a book $B$ near some
point, we have $\delta L=\delta B$ near that point.
\end{defin} Corner loci are connected with tropical and toric geometry by the following well known facts:
\begin{utver} \label{lcorner0} 1) Corner loci, and only they, are tropical hypersurfaces.
\newline 2) The isomorphism $I_{\Q}$ sends every
conewise linear function to its corner locus.
\end{utver}
For instance, the corner locus $(P,p)$ of the support function of
an integer polytope $A$ admits the following simple description:
the set $P$ contains all external normal covectors to the edges
of $A$, and the value of $p$ at such a covector equals the integer
length of the corresponding edge. In this case, $A$ is called the
\textit{Newton polytope} of the tropical hypersurface $(P,p)$,
and the following tropical version of the Kouchnirenko-Bernstein
theorem is well known (note the absence of assumptions of general
position):
\begin{theor}[Tropical Bernstein theorem] \label{tropb} The intersection number of $n$ tropical hypersurfaces in
$\R^n$ equals the mixed volume of their Newton polytopes, i.e. we
have
$$\delta A_1(\cdot)\cdot\ldots\cdot\delta A_n(\cdot)=(\{0\},A_1\cdot\ldots\cdot A_n).$$
\end{theor}
\begin{exa} The support function of a triangle and its corner locus
are shown on the left of Picture 2. Thus,
the pair of triangles on Picture 1 are the Newton polygons of the
tropical curves on the right of Picture 2, so the mixed area of
the triangles equals the intersection number of the curves.
\end{exa}
 
\textbf{Proof of Theorem \ref{prodsup0} for rational polytopes.}
 
The isomorphism $I_{\Q}$ maps a conewise polynomial $F$ of degree
$n$ to a 0-dimensional tropical variety $(\{0\},c_F)$, where
$0\in\R^n$ is the origin and $c_F$ is a real number, depending on
$F$. We prove that the map, sending every conewise polynomial $F$
to the constant $c_F$, is the desired function $D$, i.e.
$$I_{\Q}\Bigl(A_1(\cdot)\cdot\ldots\cdot A_n(\cdot)\Bigr)=(\{0\},A_1\cdot\ldots\cdot
A_n).\eqno{(*\,)}$$
 
For this, we firstly note that
$$I_{\Q}\Bigl(A_1(\cdot)\cdot\ldots\cdot A_n(\cdot)\Bigr)=I_{\Q}\Bigl(A_1(\cdot)\Bigr)\cdot\ldots\cdot
I_{\Q}\Bigl(A_n(\cdot)\Bigr),$$ for every collection of integer
polytopes $A_1,\ldots,A_n$, because $I_{\Q}$ is a ring
isomorphism. Secondly, by Proposition \ref{lcorner0}(2) we have
$$I_{\Q}\Bigl(A_i(\cdot)\Bigr)=\delta A_i(\cdot).$$
 
The two latter equalities together with Theorem \ref{tropb} imply
the desired equality $(*\,)$.
 
\section{Tropical varieties with polynomial weights.}
 
It turns out that $I_{\Q}$ acts on a
conewise polynomial of arbitrary degree $d$ as the $d$-th degree
of a certain corner locus operator, generalizing Definition
\ref{corner0} (see Definition \ref{defdelta} below), in the same way
as it is shown above for $d=1$. To make this precise and applicable to
non-rational polytopes and cones, we need the
notion of a tropical variety with polynomial weights, which may
be of independent interest. We introduce this notion here, and
apply it to the study of the isomorphism $I_{\Q}$ in the next
section.
 
\textbf{Weighted fans.}
 
A \textit{convex polyhedral cone} in an $m$-dimensional vector
space $M$ is an intersection of its subspace and finitely many
open half-spaces. A union $C$ of finitely many convex polyhedral
cones in $M$ is called a \textit{smooth cone} of codimension $k$,
if every its point $x$ has a neighborhood, where $C$ coincides
with an $(m-k)$-dimensional plane. This plane is denoted by $T_xC\subset M$,
and its orthogonal complement is denoted by $N_xC\subset M^*$.
 
\begin{defin}\label{defprefan}
A \textit{weighted pre-fan} of codimension $k$ in $M$ is a
pair $(P,p)$, such that the \textit{support set} $P$ is a smooth
cone of codimension $k$, and the \textit{weight} $p$ is a function
that sends every point $x\in P$, endowed with a \textit{coorientation} $\alpha\in\{$orientations of
$N_xP\}$, to a $k$-form $p(x,\alpha)\in\wedge^k M^*$, such that
\newline 1) for every linear function $l:M\to\R$, 
vanishing on $T_xP$, we have 
$p(x,\alpha)\wedge dl=0$, 
\newline 2)  $p(y,\alpha)$ is odd as a function of $\alpha$, i.e. $p(x,\alpha)+p(x,-\alpha)=0$, and
\newline 3) $p(y,\alpha)$ is a polynomial as a function of $y$ in a neighborhood of $x$.
\end{defin}
\begin{rem}By convention, the set of orientations of the 0-dimensional space is $\{-1,+1\}$.
\end{rem}
\begin{exa} Let $x_1,\ldots,x_m$ be the standard coordinates in $\R^m$, let $P$ be the set $\{x_1=\ldots=x_k=0,\; x_{k+1}>0\}$, and choose the standard coorientation $\alpha=dx_1\wedge\ldots\wedge dx_k$ on $P$. If $(P,p)$ is a weighted pre-fan, then its weight can be written as $p(x,\pm \alpha)=\pm f(x_{k+1},\ldots,x_m)dx_1\wedge\ldots\wedge dx_k$, where $f$ is a polynomial.
\end{exa} 
\begin{defin} \label{defsum}
For weighted pre-fans $(P,p)$ and $(Q,q)$ of codimension $k$
in $M$, we define the sum $(P,p)+(Q,q)$ as the pre-fan $(R_1\sqcup
R_2\sqcup R,r)$, where
 
\vspace{0.3cm} $R_1=P\setminus\overline{Q}$, and $r=p$ on $R_1$;
 
\vspace{0.3cm} $R_2=Q\setminus\overline{P}$, and $r=q$ on $R_2$;
 
\vspace{0.3cm} $R=\{x\in P\cap Q\; |\; N_xP=N_xQ\}$, and $r=p+q$
on $R$.\end{defin}
\begin{defin} A \textit{weighted fan} of codimension $k$ in
$M$ is an equivalence class of  weighted pre-fans of
codimension $k$ with respect to the following equivalence
relation: $$(P,p)\sim(Q,q)\; \Leftrightarrow\;
(P,p)+(R,0)=(Q,q)+(R,0)\; \mbox{for some } R\subset M.$$ 
\end{defin} 
 
\textsc{Example.} A 0-dimensional weighted fan in $M$ is a pair
$(\{0\},p)$, where $p$ is a pseudo-volume form 
on $M$.
 
\textsc{Example.} A weighted fan of codimension 0 in $M$ is
represented by a pair $(P,p)$, where $P\subset M$ is a union
of open polyhedral cones, and $p:P\to\R$ is locally polynomial.
 
\textbf{Tropical varieties.}
 
For a weighted fan $(P,p)$ of codimension $k$ in $M$, it is convenient to define the restriction of the weight $p$ to the boundary of a subset of $P$ as follows.
Consider a convex codimension $k$ cone $Q$ with a facet $R$ (which is a face of maximal dimension),
and pick any point $y\in\overline{R}$, in whose small neighborhood
$P$ contains $Q$. 
Every coorientation $\alpha$ on $Q$ induces the boundary coorientation 
$\beta$ on $R$, and the limit of $p(x,\alpha)$ as $x\in Q$ tends to 
$y$ is denoted by $\partial^Q_R p(y,\beta)$. 

A point $x\in M$ outside a smooth cone $P$ of codimension $k$ is 
said to be in its \textit{stable boundary}, if, in a small neighborhood
of $x$, the set $P$ coincides with a disjoint union of finitely many
codimension $k$ half-subspaces with the common boundary subspace, containing $x$. We denote the boundary subspace 
by $\partial P_x$,
the set of the half-subspaces by $\mathcal{P}_x$ (so that $P$ coincides with 
$\cup_{Q\in\mathcal{P}_x} Q$
in a small neighborhood of $x$), and the stable boundary of $P$ by $\partial P$. 
The stable boundary is a smooth cone of codimension $k+1$.

\begin{exa} In the setting of the example, preceding Definition \ref{defsum}, we have $\partial P=\{x_1=\ldots=x_{k+1}=0\}$, the boundary coorientation of $\partial P$,  corresponding to $\alpha$, is $\beta=dx_{k+1}\wedge dx_1\wedge\ldots\wedge dx_{k}$, and $\partial^P_{\partial P}p(x,\pm \beta)=\pm f(0,x_{k+2},\ldots,x_m)dx_1\wedge\ldots\wedge dx_k$.
\end{exa} 

\begin{defin} \label{deftrop} A weighted codimension $k$ fan $(P,p)$ in $M$
is called a \textit{polynomially  weighted tropical variety of codimension $k$}, if,
for every point $x\in\partial P$ and coorientation $\beta\in\{$orientations of $N_x\partial P\}$, we have $$\sum_{Q\in\mathcal{P}_x}\partial^Q_{\partial P_x} p(x,\beta)=0.$$
\end{defin}
The space of codimension $k$ tropical varieties in $M$, whose weights are locally homogeneous polynomials of degree $d$,
is denoted by $\mathcal{K}_k^d(M)$. It is a vector space with respect to summation of Definition \ref{defsum}
and multiplication $c\cdot(P,p)=(P,c\cdot p)$.
\begin{exa} Let $P$ be a union of finitely many rays $l_i$ in $\R^2$, pick a linear function $u_i:\R^2\to\R$, vanishing on $l_i$ and defining its counterclockwise coorientation, and let $(P,p)$ be a weighted fan. Then we have $p(x,du_i)=f_i(x)du_i$ for $x\in l_i$, where $f_i$ is a polynomial function on $l_i$. In this case, $(P,p)$ is a tropical variety if $\sum_i f_i(0)u_i=0$.
\end{exa}

\textbf{Corner loci.}

\begin{lemma}\label{ldelta}
1) For a weighted fan $(P,p)$ of codimension $k$ 
and its point $x\in P$ with a coorientation $\alpha$, there exists 
a unique $(k+1)$-form $\delta p(x,\alpha)$, satisfying
the equality $$\delta p(x,\alpha)\;(v_0\wedge\ldots\wedge v_k)\;=\;\partial_{v_0}p(x,\alpha)\;(v_1\wedge\ldots\wedge v_k)$$ 
for every collection of vectors $v_0\in T_xP$ and $v_1,\ldots,v_k$ in $M$
(here $\partial_{v_0}p$ is the derivative of the function $p(\cdot,\alpha)$ along the vector $v_0$).\newline
2) Consider a convex polyhedral cone $Q$ of codimension $k$, its facet $R$, and a linear function $l:M\to\R$,
vanishing on $R$. Then we have $$\delta\partial^Q_R(dl\wedge p)=-\partial^Q_R(dl\wedge\delta p).$$
\end{lemma}
Each of these statements follows from Condition 1 of Definition \ref{defprefan}. 
We omit the proof, because both implications
become linear algebraic tautologies when written in coordinates.
\begin{exa} In the setting of the example, preceding Definition \ref{deftrop}, let $l$ be equal to $x_{k+1}$, and let $(P,p)$ be a weighted fan. Then we have $$\delta p(x,\alpha)=df(x_{k+1},\ldots,x_m)\wedge dx_1\wedge\ldots\wedge dx_k,$$ $$\partial^P_{\partial P}(dl\wedge\delta p)(x,\beta)=dx_{k+1}\wedge df(0,x_{k+2},\ldots,x_m)\wedge dx_1\wedge\ldots\wedge dx_k, \mbox{ and }$$ $$\delta\partial^P_{\partial P}(dl\wedge p)(x,\beta)=df(0,x_{k+2},\ldots,x_m)\wedge dx_{k+1}\wedge dx_1\wedge\ldots\wedge dx_k.$$ Essentially, Lemma
\ref{ldelta}(1) states that the first of these expressions makes sense as a $(k+1)$-form, and Lemma
\ref{ldelta}(2) states that the two latter expressions are equal up to the sign.
\end{exa}

\begin{defin}\label{defdelta}
For a polynomially weighted tropical variety $(P,p)$, define
the form $r(x,\alpha)$ for every point $x\in\partial P$ with a coorientation $\beta\in\{$orientations of $N_x\partial P\}$
as $$r(x,\beta)=\sum_{Q\in\mathcal{P}_x}\partial^Q_{\partial P_x}\delta p(x,\beta).$$
The weighted fan $(\partial P,r)$ is called the \textit{corner locus} of $(P,p)$ and is denoted
by $\delta(P,p)$.
\end{defin}
\begin{exa} In the setting of the example after Definition \ref{deftrop}, the corner locus $\delta(P,p)$ is the point $\{0\}$ endowed with the weight $\sum_i df_i(0)\wedge du_i$. Note that $df_i(0)\wedge du_i$ makes sense, if $df_i$ is defined on the ray $l_i$, and $u_i$ vanishes on this ray.
\end{exa}
\begin{theor} The corner locus of a tropical variety is a tropical variety.
\end{theor}
\textsc{Proof.} Let $(P,p)$ be a polynomially weighted tropical variety. In order to prove Condition 1 of Definition \ref{defprefan} for $\delta(P,p)$ at a point
$x\in\partial P$, we should prove that $dl\wedge\sum_{Q\in\mathcal{P}_x}\partial^Q_{\partial P_x}\delta p(x,\alpha)=0$ for every linear function $l$, vanishing on $T_x\partial P$.
By Lemma \ref{ldelta}(2), we can rewrite this equality as $\delta(dl\wedge\sum_{Q\in\mathcal{P}_x}\partial^Q_{\partial P_x}p)=0$.
The latter equality follows from $\sum_{Q\in\mathcal{P}_x}\partial^Q_{\partial P_x}p=0$, which is the assumption of Definition \ref{deftrop} for the tropical variety $(P,p)$.

In order to prove the assumption of Definition \ref{deftrop} for $\delta(P,p)$ at a point $x\in\partial\partial P$, it is convenient
to choose a representative weighted pre-fan $(P,p)$ of the given tropical variety, such that $P$ is the preimage
of a two-dimensional fan under a surjection $M\to N$. In more detail, the following takes place in a small 
neighborhood of $x$:
\newline
1) The set $\partial\partial P$ coincides with a subspace $R\subset M$ of codimension $k+2$,
\newline
2) The set $\partial P$ coincides with a disjoint union of finitely many half-subspaces $Q_i$,
whose common boundary is $R$,\newline
3) The set $P$ coincides with a disjoint union of finitely many convex polyhedral cones $P_j$, such that
every $P_j$ has two facets, and these facets equal $Q_{j'}$ and $Q_{j''}$ for some $j'$ and $j''$.

In this notation, we should prove the equality $\sum \partial^{Q_i}_R\partial^{P_j}_{Q_i}\delta p=0$,
where the sum is taken over all pairs $(i,j)$ such that $Q_i$ is a facet of $P_j$. To prove this
equality, sum up the tautological equalities 
$$\partial^{Q_{j'}}_R\partial^{P_j}_{Q_{j'}}+\partial^{Q_{j''}}_R\partial^{P_j}_{Q_{j''}}=0$$
over all $j$. $\quad\Box$

\textbf{Remarks}.

\begin{rem} If $M$ is endowed with a metric or with a lattice, then, identifying
vectors with covectors and pseudovolumes with scalars, weights of weighted fans
can be considered number-valued, rather than form-valued. When writing weights
as number-valued functions in subsequent examples, we
always imply that such identification is done. \end{rem}
\begin{rem} Although we only admit piecewise polynomial weights
for weighted fans, everything will work fine with piecewise smooth
weights as well. One example of where piecewise smooth weights
are relevant is kindly provided by D. Siersma. If $F(x)$ is the
distance from a point $x\in\R^n$ to a finite set $A\subset\R^n$,
then the function $F:\R^n\to\R$ is piecewise smooth, and its
$k$-th corner locus $\delta^kF$ is a well defined tropical variety
$(P,p)$. One can easily verify that $P$ is the codimension $k$
skeleton of the Voronoi diagram of $A$, and critical points of
$p$ coincide with those of the distance function $F$ contained in
$P$.
\end{rem} Many assertions in what follows are straightforward generalizations to the case of
polynomial weights of what is known about conventional tropical
varieties with constant weights. Since the proof of such
assertions repeats the case of constant weights word by word, we
omit it and refer the reader to canonical papers like
\cite{fs}, \cite{kaz} or \cite{m} for details. The only sources of new
information are the assertions about the corner locus
differential $\delta$.
\begin{lemma}\label{ldif}
1) A weighted fan $(P,p)$ of codimension 0 is a tropical variety with polynomial
weights in $M$, if and only if $P$ is a union of open (codimension 0) polyhedral cones, 
and the function $p:P\to\R$ is the restriction of a
continuous conewise-polynomial function $M\to\R$, vanishing outside
the closure of $P$.
 
2) The map
$\delta:\mathcal{K}_0^1(M)\to\mathcal{K}_1^0(M)$
is surjective with kernel $\{(M,l)\; |\; l$ is a linear function on $M\}$.
\end{lemma}
Part 2 is a new formulation of Proposition \ref{lcorner0}.
 
\textsc{Proof of Part 1.} Continuity of $p$ at points of
$\partial P$ is equivalent to the assumption of Definition \ref{deftrop} for $(P,p)$. 
Continuity at other points follows
from a toy version of the Riemann removable singularity theorem:
if a real piecewise-polynomial function is continuous outside of
a set of codimension 2, then it is continuous everywhere.
$\quad\Box$
  
\textbf{Products.}

\begin{defin} Let $(P,p)$ and $(Q,q)$ be two weighted fans,
such that the planes $T_xP$ and $T_xQ$ are transversal at some point $x\in P\cap Q$.
Then orientations $\alpha$ and $\beta$ on the spaces $N_xP$ and $N_xQ$
induce the orientation $\alpha\wedge\beta$ on $N_x(P\cap Q)=N_xP\oplus N_xQ$,
and we define the \textit{exterior product} $p\wedge q$ of the weights $p$ and $q$
at the point $x$ by the equality $$(p\wedge q)(x,\alpha\wedge\beta)=p(x,\alpha)\wedge q(x,\beta).$$
\end{defin}
 
The Cartesian product of weighted fans $(P,p)$ in $M$ and $(Q,q)$
in $N$ is the weighted fan $(P\times Q,p\wedge q)\in M\oplus N$.
It is denoted by $(P,p)\times(Q,q)$.
\begin{lemma} $ $\newline 1) If $F$ and $G$ are polynomially weighted tropical
varieties, then so is $F\times G$. \newline 2) In this case, we
have the Leibnitz rule $\delta(F\times G)=(\delta F)\times
G+F\times(\delta G)$.
\end{lemma}
We omit the proof, because both statements follow by definition.
 
A pair of smooth cones in $M$ is said to be \textit{bookwise}, if
they are preimages of smooth cones of complementary dimension in
a vector space $N$ under a projection $M\to N$, and their union
is not contained in a hypersurface. A point $x\in
\overline{P}\cap\overline{Q}$ is said to be in the \textit{stable
intersection} $P\cap_sQ$ of smooth cones $P$ and $Q$ in $M$, if,
in a small neighborhood of $x$, the pair $(P, Q)$ coincides with
a bookwise pair of cones. $P\cap_sQ$ is a smooth cone of codimension $\codim P+\codim Q$. \par
 
 In a neighborhood of $x\in P\cap_sQ$, the smooth cones
$P$ and $Q$ split into the union of their connected components
$\sqcup_iP_i$ and $\sqcup_jQ_j$ respectively. Pick a small
(relatively to the radius of the neighborhood) vector
$\varepsilon\in M$ in general position with respect to $P$ and
$Q$, and define $\varepsilon_{i,j}= \left\{
\begin{array}{ll}
                           1 & \mbox{if } P_i+\varepsilon \mbox{ intersects } Q_j \\
                           0 & \mbox{otherwise} \\
                             \end{array}
                        \right.$
(the assumption of general position is that the intersections
$(P_i+\varepsilon)\cap Q_j$ are transversal, and $\overline{P}\cap\overline{\partial
Q}=\overline{\partial P}\cap \overline{Q}=\varnothing$ in the neighborhood of $x$).
 
If $P$ and $Q$ are the support sets of weighted fans $(P,p)$ and
$(Q,q)$, then denote the limits of $p(y)$ and $q(z)$, as $y\in
P_i$ and $z\in Q_j$ tend to $x$, by $p_i$ and $q_j$ respectively.
Denote the sum $\sum_{i,j}\varepsilon_{i,j}\cdot p_i\wedge q_j$
by $s(x)$ for every $x\in P\cap_sQ$. 
\begin{defin} \label{defmult} The weighted fan $(P\cap_sQ, s)$ is called the
\textit{intersection product} of the weighted fans $(P,p)$ and
$(Q,q)$, and is denoted by $(P,p)\cdot(Q,q)$.
\end{defin}
\begin{lemma}
\label{lmult} 1) If $F$ and $G$ are polynomially weighted tropical
varieties, then so is $F\cdot G$, and its definition does not
depend on the choice of $\varepsilon$.
 
2) 
Intersection product is associative.
\end{lemma} We omit the proof as it repeats the one for tropical
varieties with constant weights.
 
\textbf{Restrictions.}

We are particularly interested in the following special case of
the intersection product. \begin{defin} Let $F$ be a polynomially
weighted tropical variety in $M$, and $L\subset M$ be a vector subspace
of codimension $d$. Choose an arbitrary constant non-zero weight $w$
such that $(L,w)$ is a tropical variety, and denote the intersection product of 
$F$ and $(L,w)$ by $(P,p)$. Then the pair $(P,p/w)$ can
be regarded as a polynomially weighted tropical variety in $L$,
does not depend on the choice of $w$, is said to be the
\textit{restriction of} $F$ to the plane $L$, and is denoted by
$F|_L$.
\end{defin}
Lemma \ref{lmult}(2) specializes to this case as follows:
\begin{lemma} \label{rr}
For any vector subspaces $K\subset L\subset M$, we have
$(F|_L)|_K=F|_K$.
\end{lemma}

\begin{theor} \label{threstr} We have $\delta(F|_L)=(\delta F)|_L$.
\end{theor}
\textsc{Proof.} If the statement is proved for $L$ being a
hyperplane, then, in general case, we can choose a complete
flag $L=L_d\subset L_{d-1}\subset\ldots\subset L_0=M$ and observe
that
$$\delta(F|_L)=\Bigl(\delta(F|_{L_d})\Bigr)\bigr|_{L_d}=\Bigl(\delta(F|_{L_{d-1}})\Bigr)\bigr|_{L_d}=
\ldots=\Bigl(\delta(F|_{L_0})\Bigr)\bigr|_{L_d}=(\delta F)|_L$$
by Lemma \ref{rr}. Thus, without loss in generality, we assume
in what follows that $L$ is a hyperplane, given by a linear equation $l=0$.

In order to prove the equality $\delta\Bigl((P,p)|_L\Bigr)=\Bigl(\delta(P,p)\Bigr)|_L$ 
near a point $x\in L\cap_s\partial P$, it is convenient
to choose a representative of the given tropical variety to be a weighted pre-fan $(P,p)$, such that $P\cap\{l>0\}$ is the preimage
of a two-dimensional fan under a surjection $M\to N$. In more detail, the following takes place in a small 
neighborhood of $x$:
\newline
1) The set $L\cap_s\partial(P\cap\{l>0\})$ coincides with a subspace $R\subset L$,
\newline
2) The set $\partial(P\cap\{l>0\})$ coincides with a disjoin union of finitely many half-subspaces $Q_i\subset M$,
whose common boundary is $R$, \newline
3) The set $P\cap\{l>0\}$ coincides with a disjoint union of finitely many convex polyhedral cones $P_j$, such that
every $P_j$ has two facets, which equal $Q_{j'}$ and $Q_{j''}$ for some $j'$ and $j''$.

In this notation, we should prove the equality 
$$\sum_{(i,j) \mbox{ \scriptsize such that }\atop Q_i\subset L\mbox{ \scriptsize is a facet of } P_j}\partial^{Q_i}_R\delta\partial^{P_j}_{Q_i}(dl\wedge p)=
\sum_{(i,j) \mbox{ \scriptsize such that }\atop Q_i\not\subset L\mbox{ \scriptsize is a facet of } P_j}\partial^{Q_i}_R\partial^{P_j}_{Q_i}(dl\wedge\delta p).$$
By Lemma  \ref{ldelta}(2), it can be rewritten as
$\sum\partial^{Q_i}_R\partial^{P_j}_{Q_i}(dl\wedge\delta p)=0$,
where the sum is taken over all pairs $(i,j)$ such that $Q_i$ is a facet of $P_j$. To prove this
equality, sum up the tautological equalities 
$$\partial^{Q_{j'}}_R\partial^{P_j}_{Q_{j'}}+\partial^{Q_{j''}}_R\partial^{P_j}_{Q_{j''}}=0$$
over all $j$. $\quad\Box$
 
\textbf{Differential ring of polynomially weighted tropical varieties.}
 
The operation of intersection product can be expressed in terms
of Cartesian product and restriction as usual:
\begin{lemma}\label{cartint} Identifying the diagonal $D$ of the
sum $M\oplus M$ with the space $M$ itself, we have $(F\times
G)|_D=F\cdot G$ for every pair of polynomially weighted tropical
varieties $F$ and $G$ in $M$.\end{lemma} We omit the proof,
because it follows by definition.
\begin{theor} \label{lleibn} If $F$ and $G$ are polynomially weighted tropical varieties
in $M$, then $\delta(F\cdot G)=\delta F\cdot G+F\cdot\delta G$.
\end{theor}
\textsc{Proof.} By Lemma \ref{cartint}, the general case can be
reduced to the case of $G=(L,c)$, where $L\subset M$ is a vector
subspace and the weight $c$ is a constant. This special case
constitutes the assertion of Theorem \ref{threstr}. $\quad\Box$
 
Let $\mathcal{K}_k^d$ be the space of all polynomially weighted
tropical varieties $(P,p)$ in the vector space $M$, such that
$\codim P=k$, and $p$ is locally a homogeneous polynomial of
degree $d$. The direct sum of the spaces $\mathcal{K}_k^d$ over
all $d\geqslant 0$ and $k=0,\ldots,m$ is denoted by $\mathcal{K}$
and is called the \textit{ring of tropical varieties with
polynomial weights}. We summarize the results of this section as
follows.
\begin{sledst} $\mathcal{K}=\bigoplus\mathcal{K}_k^d$ is a bigraded differential ring with the multiplication
$$\cdot\, :\; \mathcal{K}_k^c\oplus\mathcal{K}_l^d\to\mathcal{K}_{k+l}^{c+d}$$
of Definition \ref{defmult} and the corner locus derivation
$$\delta\, :\; \mathcal{K}_k^d\to\mathcal{K}_{k+1}^{d-1}$$
of Definition \ref{defdelta}.
\end{sledst}
 
\section{The isomorphisms.}
 
Denote the subring $\bigoplus_d\mathcal{K}_0^d$ of $\mathcal{K}$
by $\mathcal{P}$, and the subring $\bigoplus_k\mathcal{K}_k^0$ by
$\mathcal{C}$. Recall that all elements of $P$ have the form $(M\setminus\Sigma,f)$, where
$M$ is the ambient vector space of dimension $m$, the function $f:M\to\R$ is continuous and conewise-polynomial, and $\Sigma$ is the set of points where $f$ is not smooth. Thus, we will always identify $\mathcal{P}$ with the ring of
continuous conewise-polynomial functions on $M$. In $\mathcal{P}$, consider the ideal
$\mathcal{L}$, generated by all linear functions on $M$. If the vector
space $M$ is endowed with an $m$-dimensional integer lattice,
then, restricting our consideration to weighted cones, whose
support sets are unions of rational polyhedral cones, we obtain
subrings $\mathcal{K}_{\Q}, \mathcal{P}_{\Q}, \mathcal{C}_{\Q}, \mathcal{L}_{\Q}$
of the rings $\mathcal{K}, \mathcal{P}, \mathcal{C}, \mathcal{L}$. Since, in the presence of the lattice,
pseudovolumes are identified with
scalars, this definition of the rings $\mathcal{C}_{\Q},
\mathcal{P}_{\Q}$ and $\mathcal{L}_{\Q}$ agrees with the one
given in Section 1.
 
We give a combinatorial (i.e. not involving geometry and topology
of toric varieties) description of the isomorphism
$I:\mathcal{P}/\mathcal{L}\to\mathcal{C}$ and its specialization
$I_{\Q}:\mathcal{P}_{\Q}/\mathcal{L}_{\Q}\to\mathcal{C}_{\Q}$, which in particular gives
a new explicit formula for the mixed volume of polytopes in terms of the product of their support functions. For
the sake of completeness, we also recall the construction of the
isomorphisms $\mathcal{H}\to\mathcal{P}_{\Q}/\mathcal{L}_{\Q}$ and
$\mathcal{H}\to\mathcal{C}_{\Q}$ (where $\mathcal{H}$ is the
direct limit of the cohomology rings of $m$-dimensional toric
varieties, as explained in Section 1). In the next section, we discuss
what happens to the isomorphism $I:\mathcal{P}/\mathcal{L}\to\mathcal{C}$,
if we replace the ambient vector space $M$ with a tropical variety.
 
\textbf{Isomorphism $\mathcal{P}/\mathcal{L}\to\mathcal{C}$.}
 
Define the map $I:\mathcal{P}\to\mathcal{C}$ on $\mathcal{K}_0^d$
as $\delta^d/d!$.
\begin{theor}\label{mainth} We have $I(\mathcal{L})=0$, and $I:\mathcal{P}/\mathcal{L}\to\mathcal{C}$ is a ring isomorphism.\end{theor}
\begin{rem} If we pick a simple fan $\Delta$, and restrict our
consideration to polynomially weighted tropical varieties, whose
support sets are unions of cones from $\Delta$, then the
statement remains valid, and the proof is the same. \end{rem}
\begin{rem} Although the linear map
$\delta^d:\mathcal{K}_0^k\to\mathcal{K}_d^{k-d}$ is surjective
for $d=k$, and the kernel of
$\delta^d:\mathcal{K}_{k-d}^d\to\mathcal{K}_k^0$ is generated by
linear functions for $d=k$, none of this remains true for other
values of $d$. For instance, introducing the standard metric
$dx^2+dy^2$ in the coordinate plane, and thus representing
weights of plane tropical curves as real-valued functions, the
restriction of the function $|x|-|y|$ to the set $\{xy=0\}$ can
be regarded as a tropical curve $F\in\mathcal{K}_1^1$, and we
have $\delta F=0$. However, $F$ cannot be represented as the
corner locus of a conewise-quadratic function, and cannot be
represented as $\sum_i l_i F_i$ for linear functions
$l_i:\R^2\to\R$ and tropical curves $F_i$ with constant weights.
(The first statement can be verified by definition, and the
second one is true because otherwise $F=\delta(\sum_i l_i
\delta^{(-1)}F_i)$, contradicting the first statement.) It would
be interesting to explicitly describe the kernel of
$\delta^d:\mathcal{K}_{k-d}^d\to\mathcal{K}_k^0$ and the image of
$\delta^d:\mathcal{K}_0^k\to\mathcal{K}_d^{k-d}$.
\end{rem}
\textsc{Proof.} Since $\delta^{d+1}(\mathcal{K}_0^d)=0$, we have
$$\delta^{k+l}(F\cdot G)\; =\; \sum_j\; C^j_{k+l}\cdot\delta^jF\cdot\delta^{k+l-j}G\; =\; C^k_{k+l}\cdot\delta^k F\cdot\delta^l G$$
for every pair of tropical varieties $F\in\mathcal{K}_0^k$ and
$G\in\mathcal{K}_0^l$, hence $I$ is indeed a ring homomorphism.
Since $\delta(M,l)=0$ for every linear function $l$, then
$I(\mathcal{L})=0$. Since the restriction of $I$ to the degree 1
is an isomorphism $\mathcal{K}_0^1\to\mathcal{K}_1^0$ by Lemma
\ref{ldif}.2, and the ring $\mathcal{C}$ is generated by
$\mathcal{K}_1^0$ (see e.g. \cite{kaz}), then the homomorphism
$I$ is surjective.\par
 
The pairing
$F, G \mapsto F\cdot{G}$ on
$\mathcal{P}/\mathcal{L}$ is perfect (see e.g.
\cite{brion2}), i.e. the image of the component $\mathcal{K}_0^m$ in the
quotient $\mathcal{P}/\mathcal{L}$ is generated by one element $L$, and every non-zero element
$F\in\mathcal{P}/\mathcal{L}$ admits an element
$G\in\mathcal{P}/\mathcal{L}$ of complementary dimension, such
that $F\cdot G = {L}\mod\mathcal{L}$. Since $I(L)$ is non-zero in $\mathcal{C}$ by surjectivity of $I$, then $I(F)\cdot
I(G)=I(L)\ne 0$, which implies that $I(F)$ is non-zero.
Thus $I$ is injective. $\Box$
 
\textbf{Proof of Proposition \ref{prodsup01}.}
 
Introducing a metric in $\R^n$ and writing $\delta^n$ explicitly by definition, we note that the weight of the zero-dimensional tropical variety
$\delta^n(f)$ for a continuous conewise-polynomial function $f:\R^n\to\R$ is exactly the sum in the statement of Proposition \ref{prodsup01}
(note that $\delta^n(f)$ is even easier to compute, because some similar terms are collected). We can thus formulate Proposition \ref{prodsup01}
as follows.
\begin{theor} \label{mainth2} We have $$\frac{\delta^n}{n!}\Bigl(A_1(\cdot)\cdot\ldots\cdot A_n(\cdot)\Bigr)=(\{0\},A_1\cdot\ldots\cdot
A_n)$$ for every collection of polytopes $A_1,\ldots,A_n$ in
$\R^n$.
\end{theor}
\textsc{Proof.} We have
$$\frac{\delta^n}{n!}\Bigl(A_1(\cdot)\cdot\ldots\cdot A_n(\cdot)\Bigr)=I\Bigl(A_1(\cdot)\cdot\ldots\cdot A_n(\cdot)\Bigr)=I\Bigl(A_1(\cdot)\Bigr)\cdot\ldots\cdot
I\Bigl(A_n(\cdot)\Bigr),$$ for any collection of polytopes
$A_1,\ldots,A_n$, because $I$ is a ring isomorphism (see Theorem
\ref{mainth}), and is defined as $\delta^n/n!$ for a homogeneous
conewise polynomial of degree $n$. For conewise-linear functions
it is defined as $\delta$, so we have
$$I\Bigl(A_i(\cdot)\Bigr)=\delta A_i(\cdot).$$
The tropical Bernstein formula is valid for arbitrary tropical
varieties with constant weights, not only for rational ones (see
e.g. \cite{kaz}):
$$\delta A_1(\cdot)\cdot\ldots\cdot
\delta A_n(\cdot)=(\{0\},A_1\cdot\ldots\cdot A_n).$$ These three
equalities imply the desired one. $\quad\Box$
 
For instance, the mixed area of the pair of triangles on Picture 1
can be counted as follows (their support functions are denoted by
$F$ and $G$):
\begin{center} \noindent\includegraphics[width=15cm]{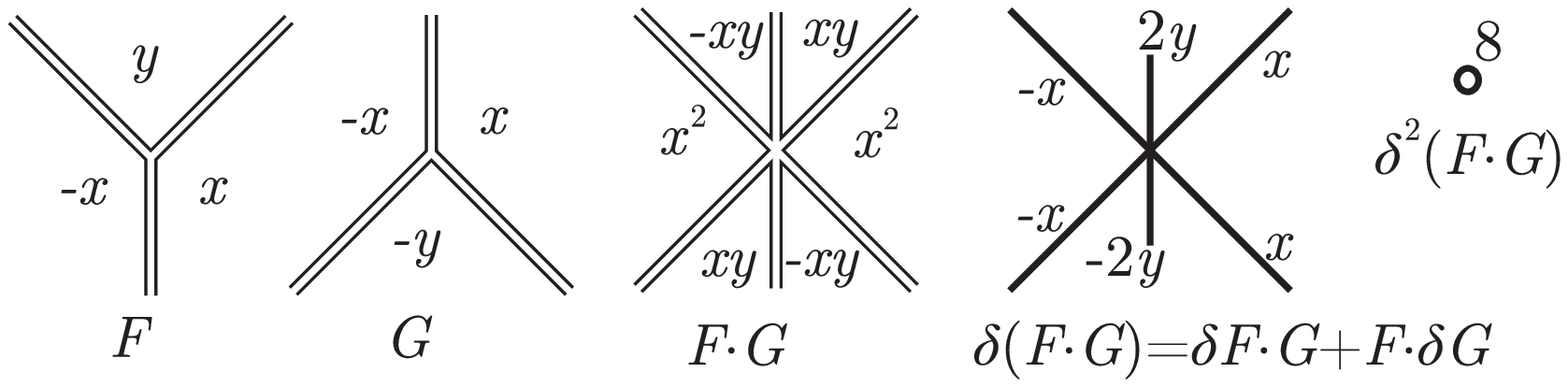}
 
\small{Picture 3.}
\end{center}
The count of the mixed area of the right pair of polygons on
Picture 1 proceeds in the same way, because the product of their
support functions is the same as for the left pair.
\begin{rem} The notion of corner loci of polynomially weighted tropical
varieties simplifies the proof of many known useful formulas for
mixed volumes. To give an example, denote the maximal face of a
polytope $A\subset\R^n$, on which a non-zero covector
$\gamma\in(\R^n)^*$ attains its maximal value $A(\gamma)$, by
$A^\gamma$, note that the $(n-1)$-dimensional mixed volume
$A_2^\gamma\cdot\ldots\cdot A_n^\gamma$ makes sense for any
polytopes $A_2,\ldots,A_n$ in the euclidean space $\R^n$, and let
$\langle\gamma\rangle$ be the ray generated by $\gamma$. Applying
the tropical Kouchnirenko-Bernstein formula to both parts of the
equality
$$\delta A_1(\cdot)\cdot\ldots\cdot \delta
A_n(\cdot)=\delta\Bigl(A_1(\cdot)\delta
A_2(\cdot)\cdot\ldots\cdot \delta A_{n}(\cdot)\Bigr),\eqno{(*)}$$
we have $\delta A_1(\cdot)\cdot\ldots\cdot \delta
A_n(\cdot)=(\{0\},A_1\cdot\ldots\cdot A_n)$ and $\delta
A_2(\cdot)\cdot\ldots\cdot \delta A_{n}(\cdot)$ is the union of
all external normal rays to the facets of $A_2+\ldots+A_n$, with
the constant weight $A_2^\gamma\cdot\ldots\cdot A_n^\gamma$
associated to every ray $\langle\gamma\rangle$. As a result, the
equality $(*)$ turns into the well known
$$A_1\cdot\ldots\cdot A_n=\sum_{|\gamma|=1}A_1(\gamma)\Bigl(A_2^\gamma\cdot\ldots\cdot
A_n^\gamma\Bigr).$$\end{rem}
 
\textbf{Isomorphisms
$\mathcal{H}\to\mathcal{P}_{\Q}/\mathcal{L}_{\Q}$ and
$\mathcal{H}\to\mathcal{C}_{\Q}$.}
 
The models $\mathcal{P}_{\Q}/\mathcal{L}_{\Q}$ and
$\mathcal{C}_{\Q}$ for the cohomology ring $\mathcal{H}$ are
Poincare dual to each other in the following sense. Pick a simple
fan $\Gamma$ in $M$, and consider a $k$-dimensional cohomological
cycle $\gamma$ in the corresponding toric variety $\T^{\Gamma}$
as an element of $\mathcal{H}$. We have the following two ways to
describe $\gamma$ explicitly. Let $\T^{C}$ be the closure of the
orbit of $\T^{\Gamma}$, corresponding to the cone $C\in\Gamma$.
The fundamental cycles of the subvarieties $\T^{C}$ over all cones
$C$ generate the homology group of $\T^{\Gamma}$, and their
Poincare duals generate the cohomology. Represent $\gamma$ as
$\sum_C\gamma_C\cdot\T^{C},\; \gamma_C\in\R$, and denote the
intersection number $\gamma\cdot\T^{C}\in\R$ by $\gamma^C$ for
every cone $C\in\Gamma$ of codimension $k$. Denote the collection of all
such cones by $\Gamma^k$. Then the cycle $\gamma$ is uniquely
determined by each of these two Poincare dual collections of
numbers
$$(\gamma_C,\; C\in\Gamma^{m-k}) \mbox{ and } (\gamma^C,\;
C\in\Gamma^{k}).$$ The image of $\gamma$ under the isomorphisms
$$I_P:\mathcal{H}\to\mathcal{P}_{\Q}/\mathcal{L}_{\Q}\quad \mbox{ and
}\quad I_C:\mathcal{H}\to\mathcal{C}_{\Q}$$ can be described in
terms of these two collections as follows.
 
For a rational subspace $L\subset\R^m$, pick a basis
$v_1,\ldots,v_l$ of the integer lattice $L\cap\Z^m$ and the
corresponding orientation $\alpha$ on $L$, and denote
$v_1\wedge\ldots\wedge v_l$ by $e(L,\alpha)$; note that
$e(L,\alpha)$ is an odd function of $\alpha$ and does not depend on the choice of $v_i$. Defining
$P=\cup_{C\in\Gamma^k} C$, and $p(x,\alpha)=\gamma^C\cdot e(N_xP,\alpha)$ for
$x\in C$, we have
$$I_C(\gamma)=(P,p).$$
 
For a simple cone $C\subset\R^m$, generated by primitive linearly
independent vectors $v_1,\ldots,v_l$, denote the polynomial
function $v^1\cdot\ldots\cdot v^l:C\to\R$ by $e(C)$, where linear
functions $v^i:C\to\R$ are dual to the vectors $v_j$ in the sense
that $v^i\cdot v_j=\delta^i_j$. Define $q(x)=\gamma_C\cdot e(C)$
for $s\in C,\; C\in\Gamma^{m-k}$, then the function $q$ on the
union $\cup_{C\in\Gamma^{m-k}} C$ admits a unique
polynomial extension of degree at most $k$ onto every cone of the
fan $\Gamma$. Gluing these extensions into a continuous
conewise-polynomial function $q:M\to\R$ of degree at most $k$,
and denoting $\cup_{C\in\Gamma^0} C$ by $Q$, we have
$$I_P(\gamma)=(Q,q).$$
 
\section{Intersection theory on tropical varieties.} 
 
We first show that the intersection theory on a smooth tropical variety
is locally induced from the ambient vector space, and then discuss
the general case. We use notation, introduced in Section 2.
 
\textbf{Intersection theory on smooth tropical varieties.}\par
A tropical variety with conewise-constant weights is considered smooth, if its support
set locally looks like a matroid fan (see e.g. \cite{fr} for the definition).
The first motivation for this terminology is to see that the tropicalization
of $V\cap\CC^n$ for an affine subspace $V\subset\C^n$ is a matroid fan.
\begin{theor} \label{smoothhh} Let the tropical variety $(P,p)$ be a matroid fan $P$ with a non-zero conewise-constant weight $p$, and suppose that $P\supset Q$ for a tropical variety $(Q,q)$.
Then $(Q,q)$ can be represented as $(P,p)\cdot V$ for some tropical
variety $V$ with conewise-polynomial weights of the same degree as $q$.
\end{theor}
\begin{rem} In this text, we restrict our attention to tropical varieties, whose support
sets consist of cones with vertices at the origin. One could also consider
``affine'' tropical varieties, whose support sets are unions of arbitrary polyhedra
of the same dimension. If we assume that $Q$ is ``affine'', then both the statement and the proof of the theorem remain valid. However, we cannot expect similar statement for ``affine'' $P$: if $P$ is the union of two parallel lines, and $Q$ is a point on one of them, then $(Q,q)=(P,p)\cdot V$ is impossible. Theorem \ref{smoothhh} is also not valid for a simplest non-smooth tropical variety (see the last example in this section). \end{rem}
The intersection theory on smooth tropical varieties, developed in \cite{fr}, \cite{all1}, \cite{sh1}, etc., is locally induced from the ambient space in the following sense: \begin{quote} The product of tropical varieties $G_1$ and $G_2$ in $(P,p)$, as defined in \cite{fr}, \cite{all1}, \cite{sh1}, equals $\widetilde G_1\cdot\widetilde G_2\cdot(P,p)$
for tropical varieties $\widetilde G_i$ such that $G_i=\widetilde G_i\cdot(P,p)$.
Such $\widetilde G_i$ always locally exist by Theorem \ref{smoothhh}.\end{quote}
In particular, the isomorphism of Theorem \ref{mainth} implies the following:
\begin{sledst} The ring of tropical varieties in a matroid fan $P$ (as constructed in \cite{fr}, \cite{all1} and \cite{sh1}) is generated by the divisors of rational functions on $P$ (in the terminology of these works).
\end{sledst}
 We recall that, for every linear map $l:M\to N$ of vector spaces, and for tropical
varieties $F$ in $M$ and $G$ in $N$, such that $\codim F \geqslant \dim\ker l$, one defines the image and the inverse image
$l_*F$ and $l^*G$, such that $l_*(F\cdot l^*G)=l_*(F)\cdot G$, and $l^*$ is a ring homomorphism (see e.g. \cite{m} or \cite{kaz}).
Let $i:M\to M\times M$ be the diagonal inclusion of the ambient
vector space $M\supset P$ in the setting of Theorem \ref{smoothhh}.
\begin{lemma}\label{lintersdiag} There exists a tropical variety $\Sigma$ in $M\oplus M$, such that, whenever $G$ is the product of tropical varieties $G_1$ and $G_2$ in $(P,p)$ in the sense of \cite{fr}, \cite{all1}, \cite{sh1}, we have $i_* G=(G_1\times G_2)\cdot\Sigma$.
\end{lemma}
\textsc{Proof.} Let $(G_1\times G_2)$ be a tropical variety $(R,r)$. In \cite{fr}, continuous conewise-linear functions $h_1,\ldots,h_k$
on the closure of $P\times P$ were constructed, such that $\delta^k (R,h_1\ldots h_kr) =i_*G$ (see Theorem 4.5 of the aforementioned work for this property). Extending the product $h_1\ldots h_k$ to a continuous conewise-polynomial function on $M\times M$, we can consider this function as
a codimension 0 tropical variety $H$ with weights of degree $k$ and take $\Sigma=\delta^kH$. $\quad\Box$\par
\textsc{Proof of Theorem \ref{smoothhh}.} Denote the tropical variety $(P,p)$ by $F$, $(Q,q)$ by $G$, and $(M,1)$ by $H$. By Lemma \ref{lintersdiag}, we have $$(F\times G)\cdot\Sigma=i_*G.$$ 
Let us now consider the diagonal inclusion
$j:M\oplus M\to(M\oplus M)\oplus(M\oplus M)$,
the projection $\pi:\oplus^3 M\to\oplus^2 M$ that sends $(b,c,d)$
to $(c,\, d-b)$, and $\pi':\oplus^4 M\to\oplus^3 M$ that sends $(a,b,c,d)$
to $(a,c,\, d-b)$, so that $\pi'=({\rm id},\pi)$. In the same way, let $i'=(i,{\rm id}): M\to\oplus^3 M$ send $a$ to $(a,a,0)$. In this notation, the aforementioned equality becomes $$(F\times G\times\Sigma)\cdot j_*(H\times H)=j_*i_*G$$ by definition of the product.
Note that $j_*(H\times H)=\pi'^*i'_*H$, thus we have $$(F\times G\times\Sigma)\cdot \pi'^*i'_*H=j_*i_*G.$$ Applying $\pi'_*$ to both sides, we have
$$\Bigl(F\times\pi_*(G\times\Sigma)\Bigr)\cdot i'_*H=i'_*G.$$
Denoting the restriction of $\pi_*(G\times\Sigma)$
to $M\oplus\{0\}\subset M\oplus M$ by $\Sigma_G$, Lemma \ref{rr} implies that
$(F\times\Sigma_G)\cdot i_*H=i_*G$, which means the desired $F\cdot\Sigma_G=G$. $\quad\Box$\par
We now axiomatize the property of $(P,p)$ that we use in the proof of Theorem \ref{smoothhh}.
 \begin{defin} A tropical variety $(P,p)$ in $M$ is said to be \textit{diagonalizable}, if it admits a tropical variety $\Sigma$ in $M\oplus M$, such that $\Bigl((P,p)\times(Q,q)\Bigr)\cdot \Sigma=i_*(Q,q)$ for every tropical variety $(Q,q)$ with $Q\subset P$.\end{defin}
 
 \begin{utver}\label{diaginters} Let the tropical variety $(P,p)$ be diagonalizable, and suppose that $P\supset Q$ for a tropical variety $(Q,q)$.
Then $(Q,q)$ can be represented as $(P,p)\cdot V$ for some tropical
variety $V$ with conewise-polynomial weights of the same degree as $q$.
\end{utver}

 
\textbf{Cohomology of tropical varieties. } \par

Intersection theory on tropical varieties (see e.g. \cite{m}, \cite{ar}, \cite{katz2}) can be formulated in our terms as follows.  
Let $F=(P,p)$ be a tropical variety with constant weights in a vector space $M$, and consider the map $m:\mathcal{K}\to\mathcal{K}$ of multiplication
by $F$, so that $m(G)=F\cdot G$ (recall that $\mathcal{K}$ is the ring of polynomially weighted tropical varieties, introduced at the end of Section 2).

\begin{defin}\label{inters111} The images $m(\mathcal{K}^0_k)$ and $m(\mathcal{K}^d_0)$ are called the \textit{homology} and
the \textit{equivariant cohomology} of $F$, and are denoted by $H_k(F)$ and $HH^d(F)$ respectively. The map $m$ brings the
ring structure of the ring $\mathcal{K}$ to the direct sums $H_\bullet(F)=\bigoplus_kH_k(F)$ and
$HH^\bullet(F)=\bigoplus_dHH^d(F)$, so that the product of $m(G_1)$ and $m(G_2)$ equals $m(G_1\cdot G_2)$.
We always consider $H_\bullet(F)$ and $HH^\bullet(F)$ as rings with respect to this ring structure,
not with respect to the one induced by the inclusions $H_\bullet(F)\subset\mathcal{K}$ and $HH^\bullet(F)\subset\mathcal{K}$.
The \textit{Poincare duality} $D_F:HH^d(F)\to H_d(F)$ is defined as $D_F(G)=\frac{1}{d!}\delta^d(G)$.
The \textit{cohomology ring} $H^\bullet(F)$ of the tropical variety $F$ is the quotient of the equivariant cohomology $HH^\bullet(F)$
by the ideal $\ker D_F$.
\end{defin}

This definition makes sense because of the following facts.

\begin{lemma} 1) $\ker D_F$ is an ideal.
\par 2) The induced map $D_F:H^\bullet(F)\to H_\bullet(F)$ is a ring isomorphism. \end{lemma} 
\textsc{Proof.} We should prove that if $D_F(g)=0$ then $D_F(g\cdot h)=0$ for every $h\in HH^c(F)$. This follows from the equality $\delta^{d+c}(gh)\cdot F=\delta^d(g)\cdot\delta^c(h)\cdot F$,
which follows from the Leibnitz rule for $\delta$ and from $\delta^{d+1}g=\delta^{c+1}h=0$.

Surjectivity and multiplicativity of $D_F$ follow from surjectivity and multiplicativity in Theorem \ref{mainth}. $\quad\Box$

\begin{exa} If $F=(M,1)$ is the vector space of dimension $m$, then $H^\bullet(F)$ and $HH^\bullet(F)$ are the direct limits
of cohomology and equivariant cohomology of $m$-dimensional toric varieties (see Section 1 for details). \end{exa}      
                 
\begin{exa} In general, the group $HH^1(F)$ is well known as the group of rational functions on $F$ (\cite{ar})
or the group of mixed Minkowski weights (\cite{katz2}), the degree 1 component of $D_F$ is the intersection map, and $H_1(F)$ is the group
of Weil divisors. Note that $H^1(F)$ is a non-trivial (in general) quotient of the group of Cartier divisors, see the second remark after
Theorem \ref{mainth} for an example. \end{exa}

\begin{exa} In our notation, the self-intersection number of the classical line $L=\{x=y,\, z=0\}$ on the tropical plane
$F=\delta\max(0,x,y,z)$ in $\R^3$ can be computed as follows. Recall that the support set $P$ of $F$ is the regular part of the singular locus
of $\max(0,x,y,z)$, and that the standard metric $x^2+y^2+z^2$ on $\R^3$ allows us to consider weights of tropical varieties as scalars. In \cite{ar}, the line $(L,1)$ is represented as $D_F(g\cdot F)$, where a continuous conewise-linear
function $g$ on $\R^3$ is uniquely defined on $P$ by the following two properties: its restriction to every connected component of $P\setminus L$
is linear, and, on the boundary of these connected components, we have
$g(1,1,1)=g(0,-1,0)=g(0,0,-1)=g(-1,-1,0)=0,\; g(1,1,0)=-1,\; g(-1,0,0)=1$. One checks by definition that $\delta(g^2\cdot F)$ is
the ray generated by $(1,1,0)$ with the linear weight $-\sqrt{2}x$ on it (this is the weight in the standard metric;
the weight in the ``integer metric'' would be $-2x$). Thus the desired self-intersection number
$L\circ L=D_F(g^2\cdot F)=\frac{1}{2}\, \delta^2(g^2\cdot F)=\frac{1}{2}\, \frac{\partial(-\sqrt{2}x)}{\partial(x/\sqrt{2})}$ equals $-1$,
which agrees with \cite{ar}. \end{exa}

Besides $H_d(F)$ and $H^d(F)$, one can consider larger groups for the tropical variety $F=(P,p)$
(they depend only on the support set $P$):
the group $\overline{H}_d(F)\supset H_d(F)$ consists of all tropical varieties with constant weights that are contained in $P$
 and have codimension $d$ in it (it is usually called the \textit{group of codimension $d$ cycles} on $F$),
the group $\overline{HH}^{d}(F)\supset HH^d(F)$ consists of all polynomially weighted tropical varieties of the form $(P,q)$ for a homogeneous
(not necessarily continuous) conewise polynomial $q$ of degree $d$ on $P$,
the Poincare dual $D_F:     \overline{HH}^{d}(F)\to\overline{H}_d(F)$ is defined by $D_F(G)=\frac{1}{d!}\delta^d(G)$,
and the group $\overline{H}^{d}(F)\supset H^d(F)$ is the quotient of $\overline{HH}^{d}(F)$ by $\ker D_F$.

These larger groups have no natural ring structure, except for the following special case.
                 \begin{defin} A 1-dimensional smooth cone is said to be \textit{regular} or \textit{regularizable}, if its rays
are the external normals to the facets of a simplex  or of a product of simplices respectively. A regular or
regularizable book is the preimage of a regular  or regularizable 1-dimensional smooth cone under
a surjection of vector spaces. A smooth cone $P$ is said to be \textit{regular} or \textit{regularizable}
in codimension 1, if it coincides with a regular or regularizable book near every point of $\partial P$.
\end{defin}
For instance, locally regularizable tropical curves are those participating in the Mikhalkin
correspondence theorem.
\begin{lemma}\label{lnorm} 1) If $(P,p)$, $(P,q)$ and $(P,r)$ are three tropical varieties with the same
support set $P$, regularizable in codimension 1, and $p$ is conewise-constant and non-zero, then
$(P,\frac{qr}{p})$ is also a tropical variety. \newline 2) If, moreover, $P$ is regular in codimension 1,
then $q/p$ is the restriction of a continuous conewise-polynomial function on the ambient space to $P$
(in particular, if $q$ is conewise-constant, then $q/p$ is constant).
\end{lemma}
\textsc{Proof.} It is enough to prove the statement for an 1-dimensional $P$. Denoting generators of its rays by $v_i$,
we rewrite the statement as follows. If the vectors $v_0,\ldots,v_m$ in $\R^m$ are the external normals to the facets of an $m$-dimensional simplex, and $\sum_i a_iv_i=\sum_i b_iv_i$, then $b_i/a_i$ does not depend on $i$. If the vectors $v_0,\ldots,v_r,\; r\geqslant m,$ in $\R^m$ are the external normals to the facets of an $m$-dimensional  product of simplices, and $\sum_i a_iv_i=\sum_i b_iv_i=\sum_i c_iv_i=0$, then $\sum_i\frac{b_ic_i}{a_i}v_i=0$. Both statements are obvious. $\quad\Box$ \par

Part 2 of this lemma shows that the following construction makes sense.
\begin{defin} Let $F=(P,p)$ be a tropical variety
with conewise-constant non-zero weight $p$ on a smooth cone $P$,
 regularizable in codimension 1. Then the \textit{large cohomology} $\overline{HH}^\bullet(F)$
is the sum $\bigoplus_d \overline{HH}^{d}(F)$ with the product 
of its elements $(P,q)$ and $(P,r)$ defined as $(P,\frac{qr}{p})$.
\end{defin}
Part 1 of the same lemma implies the equality of rings $\overline{HH}^\bullet(F)={HH}^\bullet(F)$ whenever
$P$ is regular in codimension 1, but not in general. Similarly, Theorem \ref{smoothhh}
 implies the equality of groups  $\overline{H}_\bullet(F)={H}_\bullet(F)$ whenever $F$ is smooth, but not in general. Therefore it would be interesting to know whether, for some tropical varieties $F$, the Poincare duality map
$D_F$ brings the ring structure from the large cohomology $\overline{HH}^\bullet(F)$ to $\overline{H}_\bullet(F)$,
while conventional cohomology $D_F\bigl({HH}^\bullet(F)\bigr)={H}_\bullet(F)\subsetneq\overline{H}_\bullet(F)$ 
is not enough for this purpose. For this, $D_F:\overline{HH}^\bullet(F)\to\overline{H}_\bullet(F)$ should be surjective, and its kernel should be an ideal.

The study of $\ker D_F$ is beyond the scope here; we only discuss pairwise difference between the groups $D_F(\overline{HH}^\bullet)$, $H_\bullet$ and $\overline{H}_\bullet$, because they are all different in general:
\begin{exa} Let $A$ be the union of two planes $xz=0$ in $\R^3$, and let $L$ be the $x$-coordinate line. Then $L\subset A$ cannot be represented
as the product of the tropical surface $(A,1)$ and another tropical surface with constant weights. However, this line $(L,1)\in\overline{H}_1(A)$
is Poincare dual to $(A,p)\in\overline{HH}^1(A)$, where $p:A\to\R$ equals $|y|/2$ for $z=0$ and equals $z$ for $z\ne 0$.
\end{exa}
This example implies that $H_1(A)\ne\overline{H}_1(A)$, although the Poincare duality $D_F:\overline{H}^\bullet(A)\to\overline{H}_\bullet(A)$
is still an isomorphism (by Theorem \ref{mainth}, applied to the planes $x=0$ and $z=0$). 
The many definitions and examples of this section were given to formulate and justify
the following package of conjectures (see also the second remark after Theorem \ref{mainth}):

\begin{conjec} The Poincare duality $D_F:\overline{HH}^\bullet(F)\to\overline{H}_\bullet(F)$ is surjective for every tropical variety $F$. Its kernel is an ideal, if $F$ is regularizable in codimension 1 and irreducible (in particular, the ring structure on $\overline{HH}^\bullet(F)$ induces a ring structure $\mathcal{A}$ on $\overline{H}_\bullet(F)$). The variety 
 $F$ is also diagonalizable, if it is regularizable in codimension 1  and irreducible (in particular, the ring structure on the ambient space induces a ring structure $\mathcal{B}$ on $\overline{H}_\bullet(F)$ by Proposition \ref{diaginters}). The ring structures
$\mathcal{A}$ and $\mathcal{B}$ coincide. \end{conjec}
 
\textbf{Acknowledgements.} Theorems \ref{prodsup0} and \ref{mainth}
were discussed and proved in the framework of the ``Algebra,
Geometry and Topology" seminar of the University of Toronto, lead
by A. Khovanskii in 2006; I am very grateful to A. Khovanskii,
who suggested Theorem \ref{prodsup0} and Proposition \ref{prodsup01}, to K. Kaveh, M. Mazin,
other participants of the seminar, G. Gusev, E. Katz, and J. Rau for helpful attention and fruitful
discussions. I want to thank N. A'Campo, G.-M. Greuel, D. Siersma,
O. Viro and other participants of Conference on Singularities,
Geometry and Topology in honour of the 60th Anniversary of Sabir
Gusein-Zade for many important remarks and
suggestions. The first version of the paper was greatly improved
with many examples and simplifications by B. Kazarnovskii.

\end{document}